\documentstyle[12pt]{article}
\newcommand{\la}{\langle}
\newcommand{\ra}{\rangle}

\newcommand{\ptl}{\partial}

\newcommand{\G}{\Gamma}

\newcommand{\vs}{\varsigma}

\newcommand{\ol}{\overline}

\newcommand{\kn}{\mbox{ker}}

\def\sF{\hbox{$\sc I\hskip -2.5pt F\!$}}

\def \Z{\hbox{$Z\hskip -5.2pt Z$}}

\def\qed{\hfill \hfill \ifhmode\unskip\nobreak\fi\ifmmode\ifinner
         \else\hskip5pt\fi\fi
 \hbox{\hskip5pt\vrule width4pt height6pt depth1.5pt\hskip 1 pt}}
\def\semiprod{\mbox{\rb{1.5pt}{$\ssc\,|$}$\hskip -3.pt\times$}}
\def\a{\alpha}
\def\b{\beta}
\def\d{\delta}

\def\g{\gamma}
\def\G{\Gamma}
\def\l{\lambda}

\def\si{\sigma}

\def\sc{\scriptstyle}
\def\ssc{\scriptscriptstyle}
\def\F{\hbox{$I\hskip -4pt F$}}\def \Z{\hbox{$Z\hskip -5.2pt Z$}}
\def\dis{\displaystyle}
\def\cl{\centerline}

\def\ol{\overline}

\def\rar{\rightarrow}
\def\Rar{\Rightarrow}

\def\Lar{\Leftarrow}

\def\Lra{\Leftrightarrow}

\def\bs{\backslash}
\def\hs{\hspace*}
\def\rb{\raisebox}

\def\vs{\vspace*}
\def\vsp{{}}%\vs{-2.5pt}}
\def\ra{\rangle}
\def\la{\langle}

\def\ff{{{\cal F}}} \def\WW{{\cal W}}

\begin{document}
\cl{\bf ISOMORPHISM CLASSES AND AUTOMORPHISM GROUPS}
\cl{{\bf OF ALGEBRAS OF WEYL TYPE$\ssc\,$}\footnote
{AMS Subject Classification - Primary: 17B20,
17B65, 17B67, 17B68. \\
\indent \hskip .3cm  This work is supported by NSF of China and a Fund from
National Education Department of China.}}
\vsp
\par
\vs{7pt}
\par
\centerline{Yucai Su$^\dag$$\ssc\,$\footnote
{This author was partially supported by  Academy of Mathematics and Systems Sciences
during his visit to this academy.}
 \,\,\,\, and \,\,\,\, Kaiming Zhao$^\ddag$}
\vsp
\par
\vs{-5pt}
\par
 $^\dag$Department of Applied Mathematics,  Shanghai Jiaotong University,
Shanghai 200030, P. R. China

$^\ddag$Institute of Mathematics,  Academy of Mathematics and systems Sciences,
Chinese Academy of Sciences,  Beijing 100080, P. R. China
\vsp
\par\
%\vs{5pt}
\par
{\bf Abstract} In one of our recent  papers, the associative and the Lie
algebras of Weyl type $ A[ D]= A\otimes\F[ D]$ were  defined and studied,
where $ A$ is a commutative associative algebra with an identity element over
a field $\F$ of any characteristic, and $\F[ D]$ is the polynomial algebra  of
a commutative derivation subalgebra $ D$ of $ A$. In the present paper,  a
class of the above associative and Lie algebras  $ A[ D]$ with $\F$ being a
field of characteristic $0$, $ D$ consisting of
 locally finite but not locally
nilpotent derivations of $ A$,  are studied. The isomorphism classes and
automorphism groups of these associative and Lie algebras are determined.
\vsp
\par\
%\vs{7pt}
\par
\par
The structure theory and the representation theory are two of the most
important topics in the theory of Lie algebras.
The four well-known series of infinite dimensional simple Lie algebras
of Cartan type have played important roles in the theory
of Lie algebras. Generalizations of the simple Lie algebras of Cartan type
over a field of characteristic zero have been obtained by Kawamoto [1],
Osborn [2], Dokovic and Zhao [3,4,5], Osborn and Zhao [6,7] and
Zhao [8]. Passman [9], Jordan [10]
studied the Lie algebras $ A D= A\otimes D$ of generalized
Witt type constructed from the pair of a commutative associative algebra
$ A$ with an identity element and its commutative derivation subalgebra
$ D$ over a field $\F$ of arbitrary  characteristic.
Passman proved that $ A D$
is simple if and only if $ A$ is $ D$-simple and $ A D$ acts
faithfully on $ A$ except when dim$ D=1$ and char$\F=2$.
Xu [11] studied some of these simple Lie algebras of generalized Witt type and
other generalized Cartan types Lie algebras,
based on the pairs of the tensor algebra of the group
algebra of an additive subgroup of ${\F}^n$ with the polynomial algebra
in several variables and the subalgebra of commuting locally finite
derivations.  Su, Xu and Zhang [12]
gave the structure spaces of the generalized simple Lie algebras of Witt type
constructed in [11]. We\footnote{
Su Y, Zhao K. Second cohomology group of generalized Witt type
Lie algebras and certain representations.
Submitted for publication.}
determined the second cohomology group and gave some
representations of the
Lie algebras of generalized Witt type which are some Lie algebras defined
by Passman, more general than those defined by
Dokovic and Zhao, and slightly more general than those defined by Xu.
\par
In one of our recent  papers [13], we defined
the associative and the Lie algebras of Weyl type
$ A[ D]= A\otimes\F[ D]$, where $ A$ is a commutative associative algebra
 with an identity element over a field $\F$ of any characteristic,
and $\F[ D]$ is the polynomial algebra  of
a commutative derivation subalgebra $ D$ of $ A$. and we
gave the necessary and sufficient conditions for them to be simple.
More precisely,   the associative algebras or
Lie algebras $ A[ D]$ (modular its center,  as  Lie algebra) are simple
if and only if $ A$ is $ D$-simple and $ A[ D]$ acts
faithfully on $ A$.
\par
In the present paper,  we study a class of associative and Lie algebras of the
above type $ A[ D]$ with $\F$ being a field of characteristic $0$, $ D$
consisting of locally finite but not locally nilpotent derivations of $ A$.
The isomorphism classes and automorphism groups of these associative and Lie
algebras are determined. The automorphism groups and derivations of these
algebras for a special case were obtained in [14] and [15].
\par
Throughout this paper,  we assume that $\F$ is a field of characteristic zero.
\par\
%\vs{7pt}
\par
{\bf 1 Definitions and preliminary results}
%\vs{5pt}
\par
In this paper we focus on the (associative and Lie)
algebras $ A[ D]= A\otimes\F[ D] $ defined in [13],  where
$ A$ is a  commutative associative algebra with an identity element 1
over $\F$ and $ D$ is a nonzero finite dimensional $\F$-vector
space of locally finite but not locally nilpotent
commuting $\F$-derivations of $ A$ such that
$ A$ is $ D$-simple. From [12] we know that the pairs $( A, D)$
satisfying those conditions are essentially those constructed as follows.
\par
Let $\ell_1,\ell_2$ be two nonnegative integers such that
$\ell=\ell_1+\ell_2>0$. For any integers $m,n$, denote $\ol{m,n}=\{m,\cdots,n\}$.
Take any nondegenerate additive subgroup $\G$ of $\F^{\ell}$, i.e.,
$\G$ contains an \F-basis of $\F^{\ell}$. Elements in $\G$ will be written
as $\a=(\a_1,\cdots,\a_\ell)$.
\vsp
Set
$$
J_1=\Z_+^{\ell_1}\times\{0\}^{\ell_2},\ \
J_2=\{0\}^{\ell_1}\times\Z_+^{\ell_2},\ \
J=\Z_+^\ell.
\vs{-6pt}\eqno(1)$$
Elements in $J$ will be
written as $\mu=(\mu_1,\cdots,\mu_\ell)$.
But sometimes elements in $J_1$ will also be written as
${\vec i}=(i_1,i_2,\cdots,i_\ell)$.
\vsp
Denote
$$
|\mu|=\sum_{p=1}^\ell\mu_p,\ \
%\mu!=\prod_{p=1}^\ell\mu_p!,\ \
({\ssc\,}^{\dis\mu}_{\dis\nu}{\ssc\,})=
\prod_{p=1}^\ell({\ssc\,}^{\dis\mu_p}_{\dis\nu_p}{\ssc\,}),\ \
i_{[p]}=(0,\cdots,\rb{8pt}{\mbox{$^{^{\sc p}}_{_{\dis i}}$}},0,\cdots,0),
\vsp\eqno(2)$$
for $\mu,\nu\in J,i\in\Z_+,p\in\ol{1,\ell}$
(note that $({}^i_j)=0$ if $i<j$). The value
$|\mu|$ is called the {\it level} of $\mu$.
Define a total ordering on $J$
\vsp
by
$$
\mu<\nu\Lra|\mu|<|\nu|\mbox{ or }|\mu|=|\nu|
\mbox{ but }\exists\,p\in\ol{1,\ell}\mbox{ with }
\mu_p<\nu_p\mbox{ and }\mu_q=\nu_q,q<p.
\vsp\eqno(3)$$
\par
Let $( A,\cdot)$ be the semi-group algebra
$\F[\G\times J_1]$ with the basis $\{x^{\a,\vec i}\,|\,(\a,\vec i)
\in\G\times J_1\},$ and subject to the algebraic
operation
$$
x^{\a,\vec i}x^{\b,\vec j}=x^{\a+\b,\vec i+\vec j},
\ \ \forall\,(\a,\vec i),(\b,\vec j)\in\G\times J_1.
\vsp\eqno(4)$$
Denote the identity element $x^{0,0}$ by 1, and denote $x^{\a}=x^{\a,\vec 0},\,x^{\vec i}=x^{0,\vec i}$ for
$\a\in\G,\,\vec i\in J_1$.
Define the linear transformations (derivations)
$\{\partial^-_1,\cdots,\partial^-_{\ell_1},\partial^+_1,
\cdots,\partial^+_{\ell},\ptl_1,\cdots,\ptl_\ell\}$ on $ A$
\vsp
by
$$
\partial^-_p(x^{\a,\vec i})=i_px^{\a,\vec i-1_{[p]}},
\ \ \partial^+_q(x^{\a,\vec i})=\a_q x^{\a,\vec i},\ \
\partial_p=\partial^-_p+\partial^+_p,\ \
\partial_r=\partial^+_r,
\vsp\eqno(5)$$
for all $p\in\ol{1,\ell_1},\,q\in\ol{1,\ell},r\in\ol{\ell_1+1,\ell}$.
The operators $\partial_p^-$ are called {\it down-grading operators}, and
$\partial^+_q$ are {\it grading operators}.
Note that $\partial_p$ is locally finite if $p\in\ol{1,\ell_1}$,
semi-simple if $p\in\ol{\ell_1+1,\ell}$. Set
$ D={\rm span}\{\ptl_p\,|\,p\in\ol{1,\ell}\}$.
Then we obtain the
pair $( A, D)$ where $ D$ is
a nonzero finite dimensional $\F$-vector space of locally finite
but not locally nilpotent
commuting $\F$-derivations of $ A$ such that $ A$ is $ D$-simple.
Clearly, $\{a\in A\,|\, D(a)=0\}=\F$.
\par
\vsp
Denote
$$
 D_1=\sum_{p=1}^{\ell_1}\F\ptl_p,\ \ \
 D_2=\sum_{q=\ell_1+1}^{\ell}\F\ptl_q.
\vsp
\eqno(6)$$
For any $\ptl=\sum_{p=1}^{\ell}a_p\ptl_p$ $\in D$ and
$\a=(\a_1,\cdots,\a_\ell)\in\G$,
\vsp
define
$$
\la\a,\ptl\ra=\la\ptl,\a\ra=\a(\ptl)=\sum_{p=1}^\ell a_p\a_p.
\vsp
\eqno(7)$$
Since $\G$ is a nondegenerate subgroup of $\F^{\ell}$, there
exist $\F$-basis $\a^{(1)},\cdots,\a^{(\ell)}\in\G$
of $\F^{\ell}$, and  the dual
basis $d_1,\cdots,d_\ell$ of
$ D$  with respect to (7)
such that $\la\a^{(p)},d_q\ra=\d_{p,q}$ for $p,q\in\ol{1,\ell}$.
\par
Throughout the paper, we shall fix
$\a^{(p)}$ and $d_p$, $p\in\ol{1,\ell}$.
\par
Denote by $\F[ D]$ the polynomial algebra of $ D$. Then $\F[ D]$ has
a basis
$$B=\{\ptl^\mu=\prod_{i=1}^\ell\ptl_i^{\mu_i}\,|\,\mu\in J\}
\mbox{ or }
B_1=\{d{\ssc\,}^\mu=\prod_{p=1}^\ell d_p^{\mu_p}\,|\,\mu\in J\}.
\vsp
\eqno(8)$$
The \F-vector
\vsp
space
$$
{\cal W}{\sc\!}={\sc\!}{\cal W}(\ell_1,\ell_2,\G)
{\sc\!}={\sc\!} A[ D]{\sc\!}={\sc\!}
 A{\sc\!}\otimes{\sc\!}\F[ D]{\sc\!}={\sc\!}{\rm span}
\{x^{\a,\vec i}\ptl^\mu\,|\,(\a,\vec i,\mu){\sc\!}\in
{\sc\!} \G{\sc\!}\times {\sc\!}J_1{\sc\!}\times{\sc\!} J\},
\vsp
\eqno(9)$$
forms an associative algebra with  the following product,
which is called {\it an associative algebra of Weyl type},
$$
u\partial^\mu\cdot v\partial^\nu
=\sum_{\l\in J}({\ssc\,}^{\dis\mu}_{\dis\l}{\ssc\,})u\ptl^\l(v)
\ptl^{\mu+\nu-\l},  \eqno(10)$$
for all $u,v\in  A,\ \mu,\nu\in J,$
where, by notation (2), there are only finite number
of nonzero terms in (10), and where,
$\ptl^\mu(u)=\prod_{p=1}^\ell\ptl^{\mu_p}_p(u)$.
The induced Lie algebra from the above associative algebra
is called a {\it Lie} {\it algebra
of Weyl type}, the bracket is
$$
[u\partial^\mu, v\partial^\nu]
=u\partial^\mu\cdot v\partial^\nu
-v\partial^\nu\cdot
u\partial^\mu,
$$
for all $u,v\in  A,\ \mu,\nu\in J.$
Then ${\cal W}$ acts naturally on $ A$ by
$a\ptl^\mu:x\mapsto a\ptl^\mu(x)$, it gives rise to an
(associative and Lie)
\vsp
homomorphism
$$\theta:{\cal W}\rar{\rm Hom}_{\sF}( A, A).
\vsp
\eqno(11)$$
Obviously $\F$ is contained
in the center of ${\cal W}$.
\par
{\bf Theorem 1.1} The associative algebra ${\cal W}$ and the
Lie algebra $({\cal W}/\F,[\cdot,\cdot])$ are
simple.
\par
{\bf Proof.}
By [13, Theorem 1.1], it suffices to prove that
$\F[ D]$ acts faithfully on $ A$.
Suppose $u=\sum_{\mu\in J}c_\mu d{\ssc\,}^\mu\in\F[ D]$ acts trivially
on $ A$, where
$c_\mu\in\F$, and
$N_0=\{\mu\in J\,|\,c_\mu\ne0\}$ is finite.
Take $\a=\sum_{q=1}^\ell n_q\a^{(q)}\in\G$,
where $n_q\in\Z_+$ are arbitrary.
By definition (5), we
\vsp
have
$$
0=u(x^\a)=\sum_{\mu\in N_0}c_\mu
\prod_{q=1}^\ell n_q^{\mu_q}x^\a.
\vsp
\eqno(12)$$
Thus the coefficient of $x^\a$ is zero.
But $n_q\in\Z_+$ are arbitrary, this
proves that all $c_\mu=0$,
so, $u=0$, i.e., $\F[ D]$ acts faithfully on $ A$.
\qed\par
For any monomial $u=x^{\a,\vec i}\ptl^\mu$, we call $\a$ the
{\it degree of $u$ with respect to $\G$},
$\vec i$  the {\it degree of $u$ with respect to $J_1$},
$\mu$  the {\it degree of $u$ with respect to $\ptl$} and
$|\mu|$ the {\it level of $u$ with respect to $\ptl$}.
\vsp
Let
$${\cal F}=\{u\in{\cal W}\,\,|\,\,ad(u)\,\,\mbox{\rm is locally finite on }
\WW\},
\vsp
$$
$${\cal N}=\{u\in{\cal W}\,\,|\,\,ad(u)\,\,\mbox{\rm is locally nilpotent on } \WW\}.
\vsp
$$
\par
Choose a total ordering on $\G$ compatible with its group structure.
For any basis $\{\ptl'_1,\cdots,\ptl'_\ell\}$ of $ D$ and
any $\a\in\G,\vec i\in J_1,\mu\in J$,
\vsp
set
$$
\matrix{
 A[ D]_\a={\rm span}\{x^{\tau,\vec m}\ptl'^{\nu}\,|\,\tau<\a\},
\vs{4pt}\hfill\cr
 A[ D]_{\vec i}={\rm span}\{x^{\tau,\vec m}\ptl'^{\nu}\,|\,\vec m<\vec i\},
\vs{4pt}\hfill\cr
 A[ D]_{|\mu|}={\rm span}\{x^{\tau,\vec m}\ptl'^{\nu}\,|\,|\nu|<|\mu|\}.
\hfill\cr
}
\vsp
\eqno(13)$$
{}From (10) it follows that
\vsp
that
$$
\matrix{
[ A[ D]_{\a}, A[ D]_\tau]\subseteq A[ D]_{\a+\tau},
\vs{4pt}\hfill\cr
[\ptl'^\mu,x^{\vec i}]\in A[ D]_{\vec i},
\vs{4pt}\hfill\cr
[ A[ D]_{|\mu|}, A[ D]_{|\nu|}]\subseteq A[ D]_{|\mu|+|\nu|-1}.
\hfill\cr}
\vsp
\eqno(14)$$
\par
The following lemma is crucial in obtaining
our main results in Sect.2.
\par
{\bf Lemma 1.2} (a) $\ff= D+ A$, (b) ${\cal N}= A$.
\par
{\bf Proof.} (a)
Suppose $u\!\in\!\ff\bs(D\!+\!A)$. Decompose $u$ according to the level
with respect to
\vsp
$d$,
$$u=\sum_{i=0}^n u_i,\ \
u_i=\sum_{\mu\in J,|\mu|=i}a_\mu d^\mu,a_\mu\in A,
\vsp
\eqno(15)$$
where $n>0$ is the highest level of terms in $u$, so $u_n\ne0$.
Note that, if $n=1$, we must have $u_n\notin D$.
\vsp
Write
$$
u_n\!=\!\sum_{(\a,\vec i)\in M_0}x^{\a,\vec i}u^{(\a,\vec i)}_n
\mbox{ with } u^{(\a,\vec i)}_n\in\F[ D].
\vsp
\eqno(16)$$
So
$M_0\!=\!\{(\a,\vec i)\!\in\!\G\!\times\! J_1\,|\,u^{(\a,\vec i)}_n\!\ne\!0\}$
is finite.
Let $\G_0=\{\a\,|\,\exists\,\vec i\in J_1,\,(\a,\vec i)\in M_0\}$.
Let $\b$ be the maximal element in $\G_0$. If $\G_0\ne\{0\}$, by
reversing the ordering if necessary, we can suppose $\b>0$.
Set $\vec j={\rm max}\{\vec i\in J_1\,|\,(\b,\vec i)\in M_0\}$.
\par
{\it Case 1:} $\b=0$.
\par
In this case we see that $u_n\in{\rm span}
\{x^{\vec i}\ptl^{\mu}\,|\,\vec i\in J_1,
\mu\in J\}$. Write
$$
u^{(\b,\vec j)}_n=\sum_{\mu\in J,|\mu|=n}c_\mu d^{\mu},\ c_\mu\in\F.
\vsp
\eqno(17)$$
Let $\l={\rm max}\{\mu\in J\,|\,|\mu|=n,c_\mu\ne0\}$. Suppose
$\l_k\ne0$ and $\l_i=0$ for all $i>k$.
By induction on $s$ it is easy to see that the highest term of ({\rm ad}$u)^sx^{\a^{(k)}}$ is
$\l_k^sc_{\l}^sx^{\a^{(k)},s\vec j}d^{s(\l-1_{[k]})}\ne0$.
Thus
$$\sum_{s=0}^\infty({\rm ad}u)^sx^{\a^{(k)}}=\infty.$$
this contradicts the fact that $u\in\ff$. So Case 1 does not occur.
\par
{\it Case 2:} $\b\not=0$.
\par
Choose a basis $\b^{(1)},\cdots,\b^{(\ell)}=\b$ of $\F^\ell$ in $\G$,
 and choose $\{\ptl'_1,\cdots,\ptl'_\ell\}$ to be the dual basis. Write
$u^{(\b,\vec j)}_n=\sum_{\mu\in J,|\mu|=n}c_\mu\ptl'^{\mu},\ c_\mu\in\F$.
If for all $\mu\in\{\mu\in J\,|\,|\mu|=n,c_\mu\ne0\}$, we have
$\mu_\ell=0$, then we use the argument in Case 1 to conclude
also a contradiction. So
there is a
$\mu\in\{\mu\in J\,|\,|\mu|=n,c_\mu\ne0\}$ with $\mu_\ell\ne0$.
Let $\l={\rm max}\{\mu\in J\,|\,|\mu|=n,\mu_\ell\ne0,c_\mu\ne0\}$.
It is clear that the highest term of $({\rm ad}u)^ix^{2\b}$ is
$$\b_\ell^ic_{\l}^i(2\l_\ell)(2\l_\ell+1)\cdots(2\l_\ell+
i-1)x^{(i+2)\b,i\vec j}\ptl'{}
^{i(\l-1_{[\ell]})}\ne0.$$
Thus
$$
\sum_{i=0}^\infty({\rm ad}u)^ix^{2\b}=\infty.
\vsp
$$
This contradicts the fact that $u\in\ff$. So Case 2 does not occur.
\par
Therefore $\ff\subset  D+ A$.
\par
For any $\ptl+u\in D+ A$, and any $xd^{\mu}$ with $x\in \F [\Gamma\times
J_1]$. Let $$ U={\rm
span}\{\ptl^m(x)(\ptl^{n_1}d^{\b_1})(u)(\ptl^{n_2}d^{\b_2})(u)\cdots
(\ptl^{n_k}d^{\b_k})(u)d^{\mu-\b_1-\b_2-\cdots-\b_k}
 \vsp \hskip 3.5cm$$ $$
\hskip 4cm \,\,|\,\, m,n_1,\cdots n_k\in\Z_+, \mu-\b_1-\b_2-\cdots-\b_k\ge 0\}.
\vsp
$$
Since the action of $D$ on $A$ is locally finite, it is easy to verify
that dim$U<\infty$ and ${\rm ad}(\ptl+u))^i(xd^{\mu})\in U$. So $\ptl+u\in\ff$,
i.e., $ D+ A\subset\ff$. Thus (a) follows.
\par
(b) For any $\ptl+u\in D+ A$ with $\ptl\ne0$, choose $\a\in\Gamma$ such that
$\ptl(\a)\ne0$. From (ad($\ptl+u))^ix^{\a}=\ptl(\a)^ix^{\a}\ne0$ we see that
$\ptl+u\notin\cal N$. On the other hand, for any $u\in A$, and any
$xd^{\mu}\in\WW$, from (14) it follows that
$({\rm ad}u)^{|\mu|+1}(xd^{\mu})=0$. Thus (b) follows.\qed
\par\
%\vs{7pt}
\par
{\bf 2 Isomorphism classes and automorphism groups}
%\vs{5pt}
\par
Now we are ready to prove the isomorphism theorem.
\par
{\bf Theorem 2.1}.
Let ${\cal W}={\cal W}(\ell_1,\ell_2,\G)$,
${\cal W}'={\cal W}(\ell'_1,\ell'_2,\G')$
be two (associative or Lie) algebras of Weyl type. Then
there exists an (associative or Lie) isomorphism
$\si:{\cal W}\cong{\cal W}'$ if and only if
there exist a group isomorphism $\tau:\G\cong\G'$ and a space linear
isomorphism $\phi: D\rar D'$ such that
$\phi( D_2)= D'_2$
\vsp
and
$$
(\ell_1,\ell_2)=(\ell'_1,\ell'_2),\ \
\la\a,\ptl\ra=\la\tau(\a),\phi(\ptl)\ra,\ \forall\,
\a\in\G,\,\ptl\in D.
\vsp
\eqno(18)$$
\par
{\bf Proof.}
We shall use the same notation but with a prime to denote elements
associated with ${\cal W}'$.
%In the following we shall use
%$ A[ D]$ to mean the associative algebra $({\cal W},\cdot)$.
\par
``$\Lar$'':  By assumption, there exists a nondegenerate
$\ell\times\ell$ matrix $G=({}^{M\ 0}_{P\ Q})\in GL_\ell(\F)$,
where $M\in GL_{\ell_1}(\F),Q\in GL_{\ell_2}(\F)$, such
\vsp
that
$$
\phi(\ptl_1,\cdots,\ptl_\ell)=(\ptl'_1,\cdots,\ptl'_\ell)G.
\vsp
\eqno(19)$$
For $\a=(\a_1,\cdots,\a_\ell)\in\G$, since $\a_p=\la\a,\ptl_p\ra
=\la\tau(\a),\phi(\ptl_p)\ra,p\in\ol{1,\ell}$, we then
\vsp
have
$$
\tau(\a)=\a'=(\a'_1,\cdots,\a'_\ell)=
(\a_1,\cdots,\a_\ell)G^{-1}=\a G^{-1}.
\vsp
\eqno(20)$$
Define a linear
\vsp
map:
$$
\matrix{
\si:\!\!\!\!&x^\a\mapsto x'{\ssc\,}^{\tau(\a)},\
\ptl_q\mapsto\phi(\ptl_q),\ \
\forall\,\a\in\G,\,q\in\ol{1,\ell},
\vs{4pt}\hfill\cr
&(x^{1_{[1]}},\cdots,x^{1_{[\ell_1]}})\mapsto
(x'{\ssc\,}^{1_{[1]}},\cdots,x'{\ssc\,}^{1_{[\ell_1]}})
(M^{{\ssc\,}\rm t})^{-1},
\hfill\cr
}
\vsp
\eqno(21)$$
where $M^{{\ssc\,}\rm t}$ is the transpose of $M$. It is straightforward
to verify
\vsp
that
$$\si(\ptl_p)(\si(x^{1_{[q]}}))=\d_{p,q}=\ptl_p(x^{1_{[q]}})
\mbox{ and }\si(\ptl_p)(\si(x^{\a}))=\si(\ptl_p(x^{\a})).
\vsp
\eqno(22)$$
Since the associative algebra $ A[ D]$
is generated by the elements appeared in (21), thus (21) uniquely
extends to an associative isomorphism
$\si: A[ D]\cong A'[ D']$. Thus they are also isomorphic as
Lie algebras.
\par
``$\Rar$'':
If ${\cal W}$ and ${\cal W}'$ are isomorphic as associative algebras,
then they must be isomorphic as Lie algebras, so we suppose they are
isomorphic as Lie algebras.
\par
Let $\ptl\in D$.
Since $\ptl$ is locally finite on $ A$, ${\rm ad}{\ptl}$ is locally
finite on ${\cal W}$, and thus ${\rm ad}{\si(\ptl)}$ is locally finite
on ${\cal W}'$, by Lemma 2.2, $\si(\ptl)\in A'+ D'$. Thus
we can write $\si(\ptl)=\phi(\ptl)+a_{\ssc\ptl}$
with $\phi(\ptl)\in D',a_{\ssc\ptl}\in A'$ such that $\phi: D\rar D'$
is a linear map. Since each $\ptl\ne0$ is not ad-locally nilpotent,
by Lemma 2.2 again, we see that $\phi$ is injective.
\par
Similarly, by exchanging position of ${\cal W}$ and ${\cal W}'$,
there exists a linear injection $\phi': D'\rar D$ such that for
all $\ptl'\in D'$, there exists $b_{\ssc\ptl'}\in A$ such that
$\si^{-1}(\ptl')=\phi'(\ptl')+b_{\ssc\ptl'}$.
\vsp
Then
$$
\ptl'=
\si(\si^{-1}(\ptl'))=\si(\phi'(\ptl')+b_{\ssc\ptl'})
=\si(\phi'(\ptl'))+\si(b_{\ssc\ptl'})
=\phi(\phi'(\ptl'))+a_{\ssc\phi'(\ptl')}+\si(b_{\ssc\ptl'}).
\vsp
\eqno(23)$$
Note that $b_{\ssc\ptl'}\in A$ is ad-locally nilpotent, and so
$\si(b_{\ssc\ptl'})$ is ad-locally nilpotent, thus by Lemma 2.2,
$\si(b_{\ssc\ptl'})\in A'$. Hence by (23), we have
$\phi(\phi'(\ptl'))=\ptl'$. Thus both $\phi$ and $\phi'$ are linear
isomorphism, and so $\ell=\ell'$.
\par
For any $\ptl\in D_2$, suppose $\phi(\ptl)\notin D'_2$, say
$\phi(\ptl)=\sum_{p=1}^{\ell'}a_p\ptl'_p$ with $a_p\ne0$ for some
$p\le\ell'_1$.
\vsp
Then
$$[\si(\ptl),x'{\ssc\,}^{1_{[p]}}]=
[\phi(\ptl)+a_{\ssc\ptl},x'{\ssc\,}^{1_{[p]}}]=a_p\ne0
\mbox{ but }({\rm ad}{\si(\ptl)})^2(x'{\ssc\,}^{1_{[p]}})=0.
\vsp
\eqno(24)$$
This contradicts that $\si(\ptl)$ is ad-semi-simple.
This proves that
$\phi( D_2)\subset D'_2$ and so $\ell_2\le\ell'_2$. Similarly,
$\ell'_2\le\ell_2$ and so $\ell_2=\ell'_2$ and $\phi( D_2)= D'_2$.
This proves the first equation
of (18), and so $J'_1=J_1,J'=J$.
\par
For any $\a\in\G$, $\ptl'=\sum_{p=1}^\ell a_p\ptl'_p\in D'$. Let
$\ptl=\phi^{-1}(\ptl')$ and write $\si^{-1}(\ptl')= \ptl+b_{\ssc\ptl'}$ for
some $b_{\ssc\ptl'}\in A$. Since $\si(x^\a)$ is ad-locally nilpotent, by Lemma
2.2, we have $\si(x^\a)\in A'$. Because the vector space on which $ A+ D$
semisimply acts is span$\{x^{\a}|\a\in \G\}$, and because
$\sigma( A+ D)= A+ D$, we have $\si(x^\a)\in$span$\{x^{\a}|\a\in \G\}$. So we
can suppose $$ \si(x^\a)= \sum_{\a'\in M_\a}c_{\a'}x'{\ssc\,} ^{\a'}, \vsp
\eqno(25)$$ where $M_\a$ is a finite subset of $\G'$ such that $0\ne
c_{\a'}\in\F$ for all $\a'\in M_\a$. Then $$ \matrix{ [\ptl',\si(x^\a)]
\!\!\!\!&= \dis\sum_{\a'\in M_\a}c_{\a'} \la\ptl',\a'\ra x{\ssc\,}^{\a'}
\vs{4pt}\hfill\cr&= \si([\ptl+b_{\ssc\ptl'},x^\a]) \vs{4pt}\hfill\cr&=
\la\ptl,\a\ra\si(x^\a) \vs{4pt}\hfill\cr&= \la\ptl,\a\ra \dis\sum_{\a'\in
M_\a}c_{\a'} x'{\ssc\,}^{\a}. \hfill\cr } \vsp \eqno(26)$$ Since
$\ptl'\in D'$ is arbitrary,  then $M_\a$ is a singleton.
 Thus there exist a unique $a'\in\G'$ such \vsp that $$
\si(x^\a)=c_{\a}x'{\ssc\,}^{\a'}, \vsp \eqno(27)$$
 where $c_{\a}=c_{\a',0}\in\F\bs\{0\}$. Thus (27) defines a \vsp
map
 $$ \tau:\G\rar\G'\mbox{ such that }\tau(\a)=\a' \mbox{ and }\tau(0)=0, \vsp
\eqno(28)$$
 such that by (26), we \vsp have $$
\la\ptl,\a\ra=\la\phi(\ptl),\tau(\a)\ra,\ \forall\,\a\in\G. \vsp \eqno(29)$$
Then for any $\a,\b\in\G,\ptl\in D$, by (29), we \vsp have $$ \matrix{
\la\phi(\ptl),\tau(\a+\b)\ra\!\!\!\!& =\la\ptl,\a+\b\ra \vs{4pt}\hfill\cr&
=\la\ptl,\a\ra+\la\ptl,\b\ra \vs{4pt}\hfill\cr&
=\la\phi(\ptl),\tau(\a)\ra+\la\phi(\ptl),\tau(\b)\ra \vs{4pt}\hfill\cr&
=\la\phi(\ptl),\tau(\a)+\tau(\b)\ra. \hfill\cr} \vsp \eqno(30)$$ Since
$\phi( D)= D'$, (30) shows \vsp that $$ \tau(\a+\b)=\tau(\a)+\tau(\b),\ \
\forall\,\a,\b\in\G, \vsp \eqno(31)$$ i.e., $\tau$ is a group homomorphism.
If $\tau(\a)=0$, then by (27), $\a=0$, i.e., $\tau$ is injective. By (27),
we see that $\tau$ must be surjective, i.e., $\tau$  is a group isomorphism.
\qed\par Next we shall determine the automorphism group ${\rm Aut}({\cal
W},\cdot)$ of the associative algebra $({\cal W},\cdot)$ and the automorphism
group ${\rm Aut}({\cal W},[\cdot,\cdot])$ of the Lie algebra $({\cal
W},[\cdot,\cdot])$. It is clear that ${\rm Aut}({\cal W},\cdot)\subset{\rm
Aut}({\cal W},[\cdot,\cdot])$. We shall simply denote the former group  by
\def\AutW{\mbox{${\rm Aut}({\cal W})$}}%
\AutW.
Set $ A_1=\F[\G]={\rm span}\{x^\a\,|\,\a\in\G\}$, the group algebra
of $\G$.
Set $\F_2^{\ell_2}=\{0\}\times\F^{\ell_2}$, the subspace of $\F^\ell$ of
dimension $\ell_2$. First, denote the linear automorphism group of
$\G$
\vsp
 by
$${\rm Aut}_L(\G)=\{G\in GL_\ell(\F)\,|\,\G G=\G\},
\mbox{ where } \G G=\{\a G\,|\,\a\in\G\}.
\vsp
\eqno(32)$$
(cf. (20).$\ssc\,$)
\vsp
Take
$$
{\rm Aut}_2(\G)=\{G=({}^{M\ 0}_{P\ Q})\in {\rm Aut}_L(\G)\,|\,\G G=\G\},
\vsp
\eqno(33)$$
where $M\in GL_{l_1}(F)$ and so on.
Let $\tau\in{\rm Aut}( A_1)$ be an automorphism of $ A_1$.
Since $\tau$ maps invertible elements to invertible elements,
\vsp
and
$$
(\bigcup_{\a\in\G}\F x^\a)\bs\{0\}=
\mbox{ the set of invertible elements in } A_1,
\vsp
\eqno(34)$$
thus $\tau$ determines a unique group automorphism
$\tau^*\in{\rm Aut}(\G)$ and a multiplicative function $f_\tau:\G\rar\F^*$,
i.e.,
$f_\tau(\a+\b)=f_\tau(\a)f_\tau(\b)$, such
that $\tau(x^\a)=f_\tau(\a)x^{\tau^*(\a)}$ for all $\a\in\G$.
So ${\rm Aut}( A_1)={\rm Aut}(\G)\otimes{\rm Hom}(\G,\F^*)$, where
${\rm Hom}(\G,\F^*)$ is the group of
multiplicative functions $f_\tau:\G\rar\F^*$.
\vsp
Set
$${\rm Aut}_2( A_1)={\rm Aut}_2(\G)\otimes{\rm Hom}(\G,\F^*).
\vsp
\eqno(35)$$
Thus elements in ${\rm Aut}_2( A_1)$ can be written as $\tau=(G,f)$,
where $G=({}^{M\ 0}_{P\ Q})$. And in this notation we have
$\tau(x^{\a})=f(\a)x^{\tau^*(\a)}=f(\a)x^{\a G^{-1}}$.
For any $\tau=(G,f)\in{\rm Aut}_2( A_1)$, define
$\phi_\tau: D\rar D$ such
\vsp
that
$$
\phi_\tau(\ptl_1,\cdots,\ptl_\ell)=(\ptl_1,\cdots,\ptl_\ell)G.
\vsp
\eqno(36)$$
\vsp
Then
$$
\la\phi_\tau(\ptl),\tau^*(\a)\ra=\la\ptl,\a\ra,\ \
\forall\,\a\in\G,\,\ptl\in D.
\vsp
\eqno(37)$$
\vsp
As in (21), the map
$$
\matrix{
\si_\tau:\!\!\!\!&x^\a\mapsto f_\tau(\a)x^{\tau^*(\a)},\
(\ptl_1,\cdots,\ptl_\ell)\mapsto(\ptl_1,\cdots,\ptl_\ell)G,\
\forall\,\a\in\G,
\vs{4pt}\hfill\cr
&(x^{1_{[1]}},\cdots,x^{1_{[\ell_1]}})\mapsto
(x^{1_{[1]}},\cdots,x^{1_{[\ell_1]}})
(M^{{\ssc\!\,}\rm t})^{-1},
\hfill\cr
}
\vsp
\eqno(38)$$
defines a unique $\si_\tau\in\AutW$.
Obviously, $\si_\tau\si_{\tau'}=\si_{\tau\tau'}$ and we obtain
an group homomorphism $\pi:{\rm Aut}_2( A_1)\rar\AutW$ such that
$\pi(\tau)=\si_{\tau}.$ Clearly $\pi$ is injective.
\par
Next, for any $u\in A$, since $u$ is ad-locally nilpotent, $\si_u={\rm
exp}({\rm ad\ssc\,}u)$ is an automorphism of ${\cal W}$. Regarding $ A$ as an
additive group, then $\pi': A\rar\AutW$ with $\pi'(u)=\si_u$ is a group
homomorphism such that $\kn\pi'=\F$. It is straightforward to verify that
$\delta\si_u\delta^{-1}=\si_{\delta(u)}$ for all $\delta\in Aut({\cal W})$  and
$u\in  A$.
\par
Third, for any $v=(v_1,\cdots,v_\ell)\in\F^\ell$,
let
$v^{(1)}=(v_1,\cdots,v_{\ell_1})\in\F^{\ell_1},
v^{(2)}=(0,\cdots,0,v_{\ell_1+1},$ $\cdots,v_\ell)\in\F_2^{\ell_2}$.
the linear map\vsp
$$
\matrix{
\si_v:\!\!\!\!&x^\a\mapsto x^{\a},\
(\ptl_1,\cdots,\ptl_\ell)\mapsto(\ptl_1,\cdots,\ptl_\ell)+v^{(2)},\
\forall\,\a\in\G,
\vs{4pt}\hfill\cr
&(x^{1_{[1]}},\cdots,x^{1_{[\ell_1]}})\mapsto
(x^{1_{[1]}},\cdots,x^{1_{[\ell_1]}})+v^{(1)},
\hfill\cr
}
\vsp
\eqno(39)$$
can be uniquely extended to be an automorphism of $\cal W$.
Thus the map
$\F^\ell\rar\AutW:v\mapsto\si_v$ is an injective group homomorphism.
\par

For any $\tau=(G,f)\in{\rm Aut}( A_1)$, where $G=({}^{M\ 0}_{P\ Q})$, and for
any $v=(v_1,\cdots,v_\ell)\in \F^\ell$, set
$$\tau[v]=(v_{\ell_1+1},\cdots,v_\ell)P\left(\matrix{ x^{1_{[1]}}\cr \cdots\cr
x^{1_{[\ell_1]}} } \right)\in  A,$$ $$\tau(v)=v({}^{(M^{{\ssc\,}\rm
t})^{-1}}_{\ \,0}\ ^0_Q) \in \F^\ell.$$ It is not  difficult to verify (only need
on generators: $ D, x^{\a}, x^{1_{[p]}}$)  that
$$\si_{\tau}^{-1}\si_{v}\si_{\tau}=\si_{\tau[v]}\si_{\tau(v)}.\eqno(40)$$ For
any $\tau,\tau'\in {\rm Aut}( A_1)$, $u,u'\in  A$, $v,v'\in \F^\ell$, we deduce
that
$$
\matrix{
(\si_{\tau}\si_{u}\si_{v})(\si_{\tau'}\si_{u'}\si_{v'})\!\!\!\!&
=\si_{\tau}\si_{u}\si_{\si_{v}\si_{\tau'}(u')}\si_{v}\si_{\tau'}\si_{v'}
\vs{4pt}\hfill\cr&
=\si_{\tau}\si_{u+\si_{v}\si_{\tau'}(u')}\si_{\tau'}\si_{\tau'}^{-1}
\si_{v}\si_{\tau'}\si_{v'}
\vs{4pt}\hfill\cr&
=\si_{\tau}\si_{u+\si_{v}\si_{\tau'}(u')}\si_{\tau'}\si_{\tau'[v]}
\si_{\tau'(v)+v'}
\vs{4pt}\hfill\cr&
=\si_{\tau\tau'}\si_{\tau'^{-1}(u)+\tau'^{-1}\si_{v}\si_{\tau'}(u')}
\si_{\tau'(v)+v'}.
\hfill\cr}
\eqno(41)$$
Thus $E={\rm Aut}_2( A_1){\rm
exp}( A/\F)\si_{\sF^\ell}$ forms a subgroup of Aut$(W)$. Note that ${\rm
exp}( A/\F)\si_{\sF^\ell}$ is a normal subgroup of $E$.

If we define a group structure on the cross set ${\rm
Aut}_2( A_1)\times( A/\F)\times\F^\ell$ \vsp by $$
(\tau',u',v')\cdot(\tau,u,v)=
(\tau'\tau,\si^{-1}_{\tau}(u')+\si_{\tau'}^{-1}\si_{v'}\si_{\tau'}(u)+\tau[v'],
\tau(v')+v), $$ then this group is isomorphic to $E$.
 It is not difficult to verify that
the following injection  is a group \vsp homomorphism $$ \pi:{\rm
Aut}_2( A_1)\semiprod(( A/\F) \times\F^\ell)\rar {\rm Aut}{\cal
}W,\,\,\,(\tau,u,v)\mapsto \sigma_\tau\sigma_u\sigma_v. $$ We identify the
image of ${\rm Aut}_2( A_1)\semiprod( A/\F\oplus\F^\ell)$ under $\pi$ with
$E$.
\par
{\bf Lemma 2.2}. {The following linear map
$\si_1$ is an automorphism of
\vsp
$({\cal W},[\cdot,\cdot])$:
$$
\si_1:x^{\a,\vec i}\ptl^\mu\mapsto
-(-\ptl)^\mu\cdot x^{\a,\vec i},\ \
\forall\,(\a,\vec i,\mu)\in\G\times J_1\times J,
\vsp
\eqno(42)$$
where $(-\ptl)^\mu=\prod_{p=1}^\ell(-\ptl_p)^{\mu_p}$ and
$(-\ptl)^\mu\cdot x^{\a,\vec i}$ is the product of $(-\ptl)^\mu$ and
$x^{\a,\vec i}$.}
\par
{\bf Proof}. To verify that (42) defines an automorphism, we need the
following combinatorial identity which is not difficult to \vsp prove $$
\sum_{\l'\le\mu}(-1)^{|\l'|}(^{\hskip .2cm \mu}_{\nu-\l'})(^
{\mu+\l'-\nu}_{\,\hskip .3cm\l'})=\delta_{\nu,0}, \ \forall\,\mu,\nu\in J. \vsp
$$ Using the above identity, we \vsp have $$ \matrix{ \dis\sum_{\l\le
\mu}(_\l^\mu)(-\ptl)^{\mu-\l}\cdot \ptl^\l ( x^{\b,\vec j}) \!\!\!\!&
\dis=\sum_{\l,\,\l'\le \mu}(_\l^\mu)(^{\mu-\l}_{\
\,\l'})(-1)^{|\l'|}\ptl^{\l+\l'} ( x^{\b,\vec j})(-\ptl)^{\mu-\l-\l'}
\vs{4pt}\hfill\cr& \dis=\sum_{\nu,\,\l'\le \mu} (-1)^{|\l'|}(^{\
\,\mu}_{\nu-\l'})(^{\mu+\l'-\nu}_{\ \ \,\,\l'})\ptl^{\nu} ( x^{\b,\vec
j})(-\ptl)^{\mu-\nu} \vs{4pt}\hfill\cr& \dis=x^{\b,\vec
j}(-\ptl)^\mu,\,\,\forall\,\,\b\in\G,\mu\in J, \vec j\in J_1. \hfill\cr} \vsp
\eqno(43)$$ Using (43) we \vsp have
$$\matrix{\dis
\si_1( [x^{\a,\vec
i}\ptl^\mu,x^{\b,\vec j}\ptl^\nu])  =\!\si_1(x^{\a,\vec i}\ptl^\mu x^{\b,\vec
j}\ptl^\nu- x^{\b,\vec j}\ptl^\nu x^{\a,\vec i}\ptl^\mu)
\vs{4pt}\hfill\cr\hs{30pt}\dis
=\!\si_1(x^{\a,\vec i}\sum_{\l\le \mu}(_\l^\mu)\ptl^\l( x^{\b,\vec j})
\ptl^{\mu+\nu-\l}- x^{\b,\vec j}\sum_{\l\le \nu}(_\l^\nu)\ptl^\l(x^{\a,\vec i})
\ptl^{\mu+\nu-\l})
\vs{4pt}\hfill\cr\hs{30pt}\dis
 =\!-\sum_{\l\le \mu}(_\l^\mu)(-\ptl)^{\mu+\nu-\l}\cdot
x^{\a,\vec i} \ptl^\l( x^{\b,\vec j})+ \sum_{\l\le
\nu}(_\l^\nu)(-\ptl)^{\mu+\nu-\l}\cdot x^{\b,\vec j}\ptl^\l (x^{\a,\vec i})
\vs{4pt}\hfill\cr\hs{30pt}\dis
 =\!-(-\ptl)^{\nu}(\sum_{\l\le \mu}(_\l^\mu)(-\ptl)^{\mu-\l}\cdot\ptl^\l (
x^{\b,\vec j}))x^{\a,\vec i}+ (-\ptl)^{\mu}(\sum_{\l\le
\nu}(_\l^\nu)(-\ptl)^{\nu-\l}\cdot \ptl^\l (x^{\a,\vec i}))x^{\b,\vec j}
\vs{4pt}\hfill\cr\hs{30pt}\dis
=\!(-\ptl)^\mu \cdot x^{\a,\vec i}(-\ptl)^\nu\cdot x^{\b,\vec j}-
 (-\ptl)^\nu \cdot x^{\b,\vec j}(-\ptl)^\mu\cdot x^{\a,\vec i}
\vs{4pt}\hfill\cr\hs{30pt}\dis
 =[-(-\ptl)^\mu \cdot x^{\a,\vec i},-(-\ptl)^\nu\cdot
x^{\b,\vec j}]
\vs{4pt}\hfill\cr\hs{30pt}\dis
= \![\si_1(x^{\a,\vec i}\ptl^\mu),\si_1(x^{\b,\vec
j}\ptl^\nu)].
\hfill\cr}$$
Thus $\si_1$ is an automorphism of $({\cal W},[\cdot,\cdot])$.
\qed
\par
Clearly $\si_1$ has order 2, i.e.,
$\si_1^2=1$, and  $\si_1$ is not an automorphism of $({\cal W},\cdot)$.
%Now let $W^*$ be the group generated by $W$ and $\si_1$.
%Then $W^*\cong W\semiprod({\Z/\Z}) $.
% In this case,
%we define $f_1:\G\rar\F^*$ by $f_1(\a)=-1$ for all $\a\in\G$. If $\ell_1>0$,
%we set $\si_1=1$ and $f_1=1$. Now let
%${\rm Hom}^*(\G,\F^*)$ be the group generated by
%$\f_1$ and ${\rm Hom}(\G,\F^*)$, and take
%${\rm Aut}_2^*( A_1)={\rm Aut}_2(\G)\otimes{\rm Hom}(\G,\F^*).$
%Then we have
\par
{\bf Theorem 2.3}. $\AutW\cong W$ and
${\rm Aut}({\cal W},[\cdot,\cdot]) \cong W\semiprod\langle \si_1\rangle $.
\par
{\bf Proof}. We first show that ${\rm Aut}({\cal W},[\cdot,\cdot]) \cong W\semiprod\langle \si_1\rangle $, i.e., ${\rm Aut}({\cal W},[\cdot,\cdot])$ is generated by $W$ and $\si_1$.

 Let $\si$ be an automorphism of
$({\cal W},[\cdot,\cdot])$.
{}From the proof of  Theorem 2.1, there exist $\tau\in{\rm Aut}(\G)$
and $\phi: D\rar D$ with $\phi=\phi_\tau$ and
$\tau\in{\rm Aut}_2( A_1)$ (by (35), we regard ${\rm Aut}(\G)$ as a
subgroup of ${\rm Aut}( A_1)\ssc\,$) such that
$\si(x^{\a})\in\F x^{\tau(\a)}$ for all $\a\in \Gamma$, $\si(\ptl)\equiv\ptl\,\,{\rm mod} \, A$
for all $\ptl\in D$.
Replacing $\si$ by $\si^{-1}_\tau\si$, we see
that for any $\ptl\in D,\a\in\G$, there exist $a_{\ssc\ptl}\in A$
and $c_\a\in\F\bs\{0\}$ such
\vsp
that
$$\si(\ptl)=\ptl+a_{\ssc\ptl}\mbox{ and }
\si(x^\a)=c_\a x^\a.
\vsp
\eqno(44)$$
For $p\in\ol{1,\ell}$, denote $w_p=a_{\ssc\ptl_p}\in A$. Then by (44),
for $p,q\in\ol{1,\ell}$, we
\vsp
have
$$
0=\si([\ptl_p,\ptl_q])=[\ptl_p+w_p,\ptl_q+w_q]
=\ptl_p(w_q)-\ptl_q(w_p).
\vsp
\eqno(45)$$
\par
{\bf Claim 1}.
By replacing $\si$ with $\si'\si$ for some $\si'\in W$, we can suppose
$w_p=0,\,p=1,\cdots,\ell$.
\par
Using induction, suppose that for $q>0$, we have proved that
$w_p=0$ for all $p<q$.
\vsp
Write
$$w_q=\sum_{(\a,\vec i)\in M_q}c^{(q)}_{\a,\vec i}x^{\a,\vec i},\
w'_q=\sum_{(\a,\vec i)\in M'_q}c^{(q)}_{\a,\vec i}x^{\a,\vec i},\
w''_q=\sum_{(\a,\vec i)\in M''_q}c^{(q)}_{\a,\vec i}x^{\a,\vec i},
\vsp
\eqno(46)$$
where $M_q=\{(\a,\vec i)\in\G\times J_1\,|\,c^{(q)}_{\a,\vec i}\ne0\}$
is finite, and
$M'_q=\{(\a,\vec i)\in M_q\,|\,\a_q\ne0\}$, $M''_q=M_q\bs M'_q$.
For any $x^{\a,\vec i}$ with $\a_q\ne0$, we have
$\a_q^{-1}\ptl_qx^{\a,\vec i}-x^{\a,\vec i}=\a_q^{-1}i_qx^{\a,\vec i-1_{[q]}}$.
Using induction on $i_q$, one sees that there exists
$u'_q\in A$ such that $\ptl_q(u'_q)=w'_q$.
By (45), $\ptl_p(w'_q)=0,p<q$, thus we can deduce
$\ptl_p(u'_q)=0$ for $p<q$.
Then by replacing $\si$ by $\si_{u'_q}\si$, we still have
$w_p=0$ for $p<q$, but then $w'_q$ vanishes. Suppose $w''_q\ne0$.
If $q\le\ell_1$, then we can choose $u''_q=\sum_{(\a,\vec i)\in M''_q}
c_{\a,\vec i}(i_q+1)^{-1}x^{\a,\vec i+1_{[q]}}$. It follows that
$\ptl_q(u''_q)=w''_q$ and so by replacing $\si$ by $\si_{u''_q}\si$,
$w''_q$ vanishes also. Thus assume that $q>\ell_1$.
If $w''_q\notin\F$, then there exists some $\ptl'\in D$ such that
$[w''_q,\ptl']\ne0$, so we
\vsp
have
$$
[\ptl_q\!+\!w''_q,\ptl']\!=\!-
\ptl'(w_q'')\ne\!0,\
[\ptl_q\!+\!w''_q,[\ptl_q\!+\!w''_q,\ptl']]\!=\!
-\!\!\sum_{(\a,\vec i)\in M''_q}c_{\a,\vec i}\a_q\ptl'(x^{\a,\vec i})\!=\!0,
\vsp
\eqno(47)$$
a contradiction with the fact that $\si(\ptl_q)=\ptl_q+w''_q$ is
ad-semi-simple. Thus $w''_q\in\F$. So, replacing $\si$ by some
$\si_v\si$, we can suppose $w''_q=0$.
This complete the proof of the inductive step, and thus the claim follows.
\par
{}From Claim 1 we see that (44) becomes
 $$\si(\ptl)=\ptl,\mbox{ and }
\si(x^\a)=c_\a x^\a,
\vsp
\eqno(44')$$
 for any $\ptl\in D,\a\in\G$.
\vsp Suppose $$ \si(x^\a\ptl)=\sum c_{\b,\vec i,\mu}x^{\b,\vec i}\ptl^\mu. \vsp
\eqno(48)$$ If $c_{\b,\vec i,\mu}\ne0$ for some $\b,\vec i,\mu$ with
$|\mu|\ge2$. Then we can choose some $\g\in\G$, such that
$[\ptl^\mu,x^\g]\notin A$. Then by applying ${\rm ad\,}\si(x^\g)$ to (48) we
obtain that $\la\ptl,\g\ra\si(x^{\a+\g})\notin A$, a contradiction with
(44). Thus $|\mu|\le1$. Then we \vsp have $$ \matrix{ \la\ptl,\g\ra
c_{\a+\g}x^{\a+\g}\!\!\!&= \la\ptl,\g\ra\si(x^{\a+\g}) \vs{4pt}\hfill\cr&
=\si([x^\a\ptl,x^\g]) \vs{4pt}\hfill\cr& \dis=\sum_{|\mu|=1}c_{\b,\vec
i,\mu}x^{\b,\vec i}[\ptl^\mu,x^\g] \vs{4pt}\hfill\cr&
\dis=\sum_{|\mu|=1}c_{\b,\vec i,\mu}\la\ptl^\mu,\g\ra x^{\b+\g,\vec i},
\hfill\cr} \vsp \eqno(49)$$ holds for all $\g\in\G.$ For any $(\b,\vec
i,\mu)$ with $c_{\b,\vec i,\mu}\ne0$, by choosing $\g$ with
$\la\ptl^\mu,\g\ra\ne0$ and comparing both sides of (49), we see that
$(\b,\vec i,\mu)=(\a,0,\mu)$. Thus we can rewrite (48) \vsp as $$
\si(x^\a\ptl)=a_{\a,\ssc\ptl}+x^\a\theta_\a(\ptl), \mbox{ where
}a_{\a,\ssc\ptl}\in A,\,\theta_\a(\ptl)\in D, \vsp \eqno(50)$$ such that
$\theta_\a: D\rar D$ is a linear injection. Applying ${\rm ad}{\ptl_p}={\rm
ad}{\si(\ptl_p)}$ to (50), we obtain that
$\ptl_p(a_{\a,\ssc\ptl})=\a_pa_{\a,\ssc\ptl}$ for all $p\in\ol{1,\ell}$. \vsp
Thus $$ a_{\a,\ssc\ptl}=c_{\a,\ssc\ptl}x^\a,\ \forall\,\a\in\G. \vsp
\eqno(51)$$ Applying ad$\si(x^\b)$ to (50), using (44), we \vsp obtain $$
\la\theta_\a(\ptl),\b\ra c_\b=\la\ptl,\b\ra c_{\a+\b},\
\forall\,\a,\b\in\G,\,\ptl\in D. \vsp \eqno(52)$$ This gives
$\la\theta_\a(\ptl),\b\ra=\la\ptl,\b\ra c_{\a+\b}c_\b^{-1}$. Replacing $\b$ by
$-\b$ in this expression, and adding the two expressions, we obtain that
$\la\ptl, (c_{\a+\b}c_\b^{-1}-c_{\a-\b}c_{-\b}^{-1})\b\ra=0$ for all
$\a,\b\in\G,\ptl\in D$. Thus we \vsp obtain $$ c_{\a+\b}=c_{\a-\b}c_\b
c_{-\b}^{-1}, \vsp \eqno(53)$$ holds for all $\b\ne0$, but if $\b=0$, (53)
holds trivially. Replacing $\a$ by $\a-\b$, we \vsp obtain $$
c_{\a}=c_{\a-2\b}c_\b c_{-\b}^{-1}=c_{\a+2\b}c_{-\b} c_{\b}^{-1}, \
\forall\,\a,\b\in\G, \vsp \eqno(54)$$ where the last equality follows from
the second by replacing $\b$ by $-\b$. Replacing $\a$ by $\g$ and multiplying
the two expressions, and then setting $\g=2\b$, we \vsp obtain $$ c_\a
c_{2\b}=c_0c_{\a+2\b},\ \ \forall\,\a,\b\in\G. \vsp \eqno(55)$$ Now in
(52), replacing $\b$ by $2\b$, using (55), we \vsp obtain $$ \la
\theta_\a(\ptl)-c_0^{-1}c_\a \ptl,\b\ra=0,\ \forall\,\a,\b\in\G,\,\ptl\in D.
\vsp \eqno(56)$$ This proved that $\theta_\a(\ptl)=c_0^{-1}c_\a\ptl$. Using
this in (52), whenever $\b\ne0$, we can choose $\ptl$ with
$\la\ptl,\b\ra\ne0$, so that we \vsp obtain $$ c_\a c_\b=c_0c_{\a+\b},\ \
\forall\,\a,\b\in\G. \vsp \eqno(57)$$ Define $f:\a\mapsto c_0^{-1}c_\a$, then
(57) shows that $f\in{\rm Hom}(\G,\F^*)\subset {\rm Aut}_2( A_1)$. Thus by
replacing $\si$ by $\si^{-1}_{f}\si$, we can suppose $c_\a=c_0$ for all
$\a\in\G$ and then (56) gives $\theta_\a(\ptl)=\ptl$. So (50) becomes $$
\si(x^\a\ptl)=c_{\a,\ssc\ptl}x^\a+x^\a\ptl,\,\,\,\si(x^\a)=c_0x^\a.
\eqno(50')$$ For any $\ptl\in D$, by taking bracket of $\si(\ptl')$ for
$\ptl'\in D$ with $\si(\ptl^2)$, we see that $\si(\ptl^2)\in\F[ D]$. By
taking bracket of $\si(x^\a)$ with it, we obtain that $\si(\ptl^2)$ must take
the \vsp form $$
\si(\ptl^2)=b_{\ssc\ptl}+b'_{\ssc\ptl}\ptl+b''_{\ssc\ptl}\ptl^2,\ \
b_{\ssc\ptl},b'_{\ssc\ptl},b''_{\ssc\ptl}\in\F. \vsp \eqno(58)$$ Again,
taking bracket of $\si(x^\a)$ with (58), we \vsp obtain $$
b''_{\ssc\ptl}=c_0^{-1},\,\,\,\, 2c_{\a,\ssc\ptl}
=c_0b'_{\ssc\ptl}+(1-c_0)\la\ptl,\a\ra \mbox{ for all $\a\in\G$ with }
\la\ptl,\a\ra\ne0. \vsp \eqno(59)$$
%This shows that $c_{\a,\ssc\ptl}$ does not depend on $\a$ if $\la\ptl,\a\ra\ne0$.
%Since $c_{\a,\ssc\ptl}$ is a linear function on $\ptl$, thus
%$c_{\a,\ssc\ptl}$ is does not depend on $\a$.
By calculating the constant term of $\si([x^{-\a}\ptl,x^\a\ptl])$, we have
$c_{-\a,\ptl}=-c_{\a,\ptl}$ if $\la\ptl,\a\ra\ne0$. Using this in (59), we
obtain that $b'_{\ssc\ptl}=0$ for all $\ptl\in D$. Taking bracket of
$\si(x^\a\ptl)$ with (58) and making use of (59), we obtain \vsp that $$
2\la\ptl,\a\ra\si(x^\a\ptl^2)+\la\ptl,\a\ra^2\si(x^\a\ptl)=
c_0^{-1}(2\la\ptl,\a\ra x^\a\ptl(\ptl+c_{\a,\ssc\ptl})+ \la\ptl,\a\ra^2
x^\a(\ptl+c_{\a,\ssc\ptl})), \vsp $$ \vsp thus $$
\si(x^\a\ptl^2)=c_0^{-1}x^\a(\ptl+c_{\a,\ssc\ptl})^2 \mbox{ for all $\a\in\G$
with } \la\ptl,\a\ra\ne0. \vsp \eqno(60)$$ Then taking bracket of
$\si(x^{-\a}\ptl)$ with (60), we obtain \vsp that $$
3\la{\sc\!}\ptl,{\sc\!}\a{\sc\!}\ra (c_0^{-1}\ptl^2 \!+\!b_{\ssc\ptl})\!-\!
\la{\sc\!}\ptl,{\sc\!}\a{\sc\!}\ra^2\ptl \!=\!
c_0^{-1}\la{\sc\!}\ptl,{\sc\!}\a{\sc\!}\ra (3\ptl^2\!+\!(4c_{\a,\ssc\ptl}\!+\!
2c_{-\a,\ssc\ptl}\!-\! \la{\sc\!}\ptl,{\sc\!}\a{\sc\!}\ra)\ptl
\!+\!c_{\a,\ssc\ptl}^2\!+\!2c_{\a,\ssc\ptl} c_{-\a,\ssc\ptl}\!-\!
\la{\sc\!}\ptl,{\sc\!}\a{\sc\!}\ra c_{-\a,\ssc\ptl}), \vsp $$ thus by (59),
we \vsp obtain $$ 12c_0b_{\ssc\ptl}=(1-c_0^2)\la\ptl,\a\ra^2 \mbox{ for all
$\a\in\G$ with } \la\ptl,\a\ra\ne0. \vsp \eqno(61)$$ This shows that
$b_{\ssc\ptl}=0$ and $c_0=\pm1$. If necessary, by replacing $\si$ by
$\si_1\si$, we can suppose $c_0=1$. Then (59) gives $c_{\a,\ssc\ptl}=0$. Thus
$$ \si(x^\a\ptl)=x^\a\ptl,\,\,\,\si(x^\a)=x^\a, \eqno(50'')$$
$$\si(\ptl)=\ptl,\,\,\,\, \si(\ptl^2)=\ptl^2.\eqno(58')$$ For any $\a\in\G$.
 By Lemma 2.2, $\si(x^{\a,1_{[p]}})\in A$ for $p\in\ol{1,\ell_1}$.
Suppose $\si(x^{\a,1_{[p]}})=$\break  $\sum_{(\b,\vec i)\in M_0}b_{\b,\vec
i}x^{\b,\vec i}$ with $M_0=\{(\b,\vec i)\,|\,b_{\b,\vec i}\ne0\}$. Taking
bracket of $\si(\ptl_q),q\in\ol{1,\ell}$ with it, we \vsp obtain $$
\sum_{(\b,\vec i)\in M_0}b_{\b,\vec i} ((\b_q-\a_q)x^{\b,\vec i}+i_qx^{\b,\vec
i-1_{[q]}})-\d_{p,q}x^\a=0. \vsp \eqno(62)$$ If some $(\b,\vec i)\in M_0$
with $\b\ne\a$, we can always choose $q$ with $\b_q\ne\a_q$ and so (62) can
not hold. Similarly, if $(\a,\vec i)\in M_0$ with $\vec i\ne 0,1_{[p]}$, using
induction on $|\vec i|$, we see that (62) can not hold. Thus
$M_0\subseteq\{(\a,0),(\a,1_{[p]})\}$, i.e., we can \vsp write $$
\si(x^{\a,1_{[p]}})=b_{\a,p}x^{\a,1_{[p]}}+b'_{\a,p}x^\a, \vsp \eqno(63)$$
and (62) then gives $b_{\a,p}=1$. Taking bracket of $\si(x^\b\ptl_q)$ with
(63), we obtain that $b'_p=b'_{\a,p}$ does not depend on $\a$. Take
$v=-(b'_1,\cdots,b'_{\ell_1})\in\F^{\ell_1}$, and by replacing $\si$ by
$\si_v\si$, we can then suppose that $b'_p=0,p\in\ol{1,\ell_1}$. Then we see
that $\si$ fixes all the elements in the set
$\{x^{\b},x^{\b}\ptl,x^{\b}\ptl^2,x^{\a,1_{[p]}}\,|\,
\b\in\G,\ptl\in D,p\in\ol{1,\ell_1}\}.$ But ${\cal W}$ is generated by this
set as a Lie algebra. Thus $\si=1$, therefore ${\rm Aut}({\cal
W},[\cdot,\cdot])\cong W\semiprod\langle \si_1\rangle $.

The result $\AutW\cong W$ follows directly from ${\rm Aut}({\cal W},[\cdot,\cdot])\cong W\semiprod\langle \si_1\rangle $.
\qed\par
Finally, we remark that Theorems 2.1, 2.3 can be generalized in following
aspect.
\par
Let $\F$ be a field of characteristic zero (not necessarily algebraically
closed). Let $I_1,I$ be any two indexing sets such that $I_1\subseteq I$.
Let $\F^{I_1}=\oplus_{p\in I_1}\F$
be the direct sum of copies of $\F$ indexed by $I_1$.
Let $\G$ be a nondegenerate subgroup of $\oplus_{p\in I}\F$,
the direct sum of copies of $\F$ indexed by $I$.
Let $J_1=\oplus_{p\in I_1}\Z_+$, the direct sum
of copies of $\Z_+$ indexed by $I_1$, and
let $J=\oplus_{p\in I}\Z_+$.
Element in $\G$ and $ J$ are written
\vsp
as
$$\a=(\a_p\,|\,p\in I)\in\G,\ \
\vec i=(i_p\,|\,p\in I)\in J,
\vsp
\eqno(64)$$
where all but a finite of $\a_p=0$ and
all but a finite of $i_p=0$ and $i_p=0$ for all $p\in I\bs I_1$.
Let $ A$ be semi-group algebra
generated by $\{x^{\a,\vec i}\,|\,\a\in\G,\vec i\in J_1\}$, and
let $ A_1={\rm span}\{x^{\a}\,|\,\a\in\G\}$.
Let $ D$ be a subspace of derivations of $ A$ spanned by $\{\ptl_p\,|\,
p\in I\}$,
\vsp
where
$$
\ptl_p(x^{\a,\vec i})=\a_px^{\a,\vec i}+i_px^{\a,\vec i-1_{[p]}},
\ \ \forall\,\a\in\G,\,\vec i\in J_1,\,p\in I.
\vsp
\eqno(65)$$
Set $ D_2$ the subspace of $ D$ spanned by $\{\ptl_q\,|\,q\in I\bs I_1\}$.
\par
Then we can define as before the (associative or Lie)
\vsp
algebra
$$
{\cal W}={\cal W}(I_1,I,\G)= A[ D]
={\rm span}\{x^{\a,\vec i}\ptl^\mu\,|\,(\a,\vec i,\mu)\in \G\times J_1\times J\}
\vsp
\eqno(66)$$
and we have
\par
{\bf Theorem 2.4}. Let ${\cal W}={\cal W}(I_1,I,\G)$, ${\cal W}'={\cal
W}(I'_1,I',\G')$ be two (associative or Lie) algebras of Weyl type. Then there
exists an (associative or Lie) isomorphism $\si:{\cal W}\cong{\cal W}'$ if and
only if there exist a group isomorphism $\tau:\G\cong\G'$ and a space linear
isomorphism $\phi: D\rar D'$ such that $\phi( D_2)= D'_2$ and
$\la\a,\ptl\ra=\la\tau(\a),\phi(\ptl)\ra,$ for all $\a\in\G,\,\ptl\in D.$
Furthermore ${\rm Aut}({\cal W}(I_1,I,\G))\cong W=\! {\rm
Aut}_2( A_1)\semiprod{\rm exp}( A/\F)\si_{\sF^I}$ and ${\rm Aut}({\cal
W}(I_1,I,\G),[\cdot,\cdot])\cong W\semiprod\langle \si_1\rangle$.\qed
\par\ \par
\vs{7pt}
{{\bf References}}
 \par
\begin{description}
\item[{[1]}] Kawamoto N.  Generalizations of Witt algebras over a field of characteristic zero.
Hiroshima Math. J., 1985, 16: 417-462

\item[{[2]}] Osborn J M.  New simple infinite-dimensional Lie algebras of
 characteristic 0. J. Alg., 1996, 185: 820-835

\item[{[3]}] Dokovic D Z,  Zhao K.  Derivations, isomorphisms, and
second cohomology of generalized Witt algebras. Trans.
Amer. Math. Soc., 1998, 350(2): 643-664

\item[{[4]}] Dokovic D Z,  Zhao K.  Generalized Cartan type $W$ Lie algebras
in characteristic zero. J. Alg., 1997, 195: 170-210

\item[{[5]}] Dokovic D Z,  Zhao K.  Derivations, isomorphisms, and
second cohomology of generalized Block  algebras.
Alg. Colloq., 1996, 3(3): 245-272

\item[{[6]}] Osborn J M,   Zhao K.  Generalized Poisson bracket and
Lie algebras
of type H in characteristic 0.  Math. Z., 1999, 230: 107-143

\item[{[7]}] Osborn J M,   Zhao K.  Generalized Cartan type K Lie algebras
in characteristic 0,  Comm. Alg., 1997, 25: 3325-3360

\item[{[8]}] Zhao K. 
Isomorphisms between generalized Cartan type W Lie algebras in
characteristic zero.  Canadian J. Math., 1998,  50: 210-224

\item[{[9]}] Passman D P.  Simple Lie algebras of Witt type. J. Alg.,
1998, 206: 682-692

\item[{[10]}] Jordan D A. On the simplicity of Lie algebras of derivations of
commutative algebras. J. Alg., 2000, 228: 580-585

\item[{[11]}] Xu X. New generalized simple Lie algebras of Cartan type over a
field with characteristic 0.  J. Alg.,  2000, 224: 23-58

\item[{[12]}] Su Y, Xu X, Zhang H.  Derivation-simple algebras
and the structures of Lie algebras of Witt type. J. Alg. in press.

\item[{[13]}] Su Y, Zhao K.  Simple algebras of Weyl type. Science in China, to appear.

\item[{[14]}] Zhao K.  Automorphisms of algebras of differential operators.
J. of Capital Normal University, 1994, 1: 1-8

\item[{[15]}] Zhao K.  Lie algebra of derivations of algebras of differential
operators. Chinese Science Bulletin,  1993,  38(10): 793-798

\end{description}
\end{document}